\documentclass[12pt]{article}
\usepackage{latexsym}
\usepackage{amsfonts}
\def\abstract{\small\quotation{\hskip-\parindent\sc Abstract.}}

\def\@footnotetext#1{\insert\footins{%
\footnotesize    \interlinepenalty\interfootnotelinepenalty
    \splittopskip\footnotesep
    \splitmaxdepth \dp\strutbox \floatingpenalty \@MM
    \hsize\columnwidth \@parboxrestore
   \edef\@currentlabel{\csname p@footnote\endcsname\@thefnmark}\@makefntext
    {\rule{\z@}{\footnotesep}\ignorespaces
      #1\strut}}}
\def\abstract{\small\quotation{\hskip-\parindent\sc Abstract.}}

\def\classification{\@ifnextchar [{\@xfootnotenext}%
   {\begingroup\let\protect\noexpand
      \xdef\@thefnmark{}\endgroup
    \@footnotetext}}

\title {}

\begin{document}

\begin{center}
{\Large \bf Applications of Lie ring methods to group theory}

\bigskip

\bigskip

  {Pavel Shumyatsky}\footnote{Research supported by FAPDF and CNPq} 

        Department of Mathematics \\
University of Brasilia\\
   70910-900 Brasilia - DF, Brazil 
   
  e-mail: pavel@mat.unb.br 
\end{center} 
\bigskip

\begin{abstract}
\noindent  This article was published in 2000. Its aim is to illustrate how a Lie-theoretic
result of Zelmanov enables one to treat various problems in group theory.
\end{abstract}

\bigskip
\bigskip

\noindent{\bf 1. Introduction}
\bigskip

Lie rings were associated with $p$-groups in the 30s in the context of the
Restricted Burnside Problem and since then Lie ring methods proved to be 
an important and very effective tool of group theory. In the last 10 years the 
sphere of use of Lie rings was amplified considerably, mainly due to Zelmanov's 
outstanding contribution.  

It has been known for some time that the following assertions are equivalent.
\bigskip

\noindent{\bf 1.1.} Let $m$ and $n$ be positive integers. Then there exists
a number $B(m,n)$ depending only on $m$ and $n$ such that the order of any
$m$-generator finite group of exponent $n$ is at most $B(m,n)$.
\bigskip

\noindent{\bf 1.2.} The class of locally nilpotent groups of exponent $n$ is 
a variety.
\bigskip

\noindent{\bf 1.3.} The class of locally finite groups of exponent $n$ is 
a variety.
\bigskip

\noindent{\bf 1.4.} Any residually finite group of exponent $n$ is locally
finite.
\bigskip

The Restricted Burnside Problem is exactly the question whether
the first of the above assertions is true. In 1956 P. Hall and G. Higman reduced
the problem to the case of prime-power exponent \cite{hahi}.
This shows that 1.1 is equivalent to 1.2. In \cite{ko} A. I. Kostrikin solved
the Restricted Burnside Problem affirmatively for $n=p$, a prime. 
His proof relies on a profound study
of Engel Lie algebras of positive characteristic. The reduction to Lie algebras
was made earlier by W. Magnus \cite{ma} and independently by I. N. Sanov \cite{sa}.
Finally, in 1989, E. Zelmanov gave a complete solution of the Restricted
Burnside Problem \cite{ze1}, \cite{ze2}. In \cite{zelm} Zelmanov deduces the 
positive solution of the Restricted Burnside Problem from the following theorem.
\bigskip

\noindent{\bf Theorem 1.5.} Let $L$ be a Lie algebra generated by
$a_1,\dots,a_m$.  Suppose that $L$ satisfies a polynomial identity and each 
commutator in $a_1,\dots,a_m$ is ad-nilpotent. Then $L$ is nilpotent.
\bigskip

The above theorem proved to be a clue for some other deep results on profinite
and residually finite groups (see \cite{zelm}, \cite{ze3}).
In the present paper we will discuss some recently discovered group-theoretic
corollaries of this theorem.
Each of the principal results described
in Sections 3-5 can be viewed as a development around the Restricted Burnside
Problem in the sense of at least one of the assertions 1.1-1.4. So, Theorem 1.5
is not only one of the main tools in the solution of the Restricted Burnside
Problem but also an important
factor extending our understanding  of the problem itself. The present paper
is an attempt to make this subject accessible to a broader audience. 

In the next section  we describe the construction
that associates a Lie algebra $L(G)$ to any group $G$ and establish basic facts
reflecting relationship between $G$ and $L(G)$.  This will be a bridge linking
the above theorem of Zelmanov with group theory. 

In Section 3 we consider residually finite groups in which all commutators\linebreak
$[x_1,x_2,\dots,x_k]$ satisfy some restrictive condition.
In particular we describe a proof of the following theorem. 
\bigskip

\noindent{\bf Theorem 3.1.}
Let $k$ be an integer, $q=p^s$ a prime-power, and $G$ a residually finite group such
that $[x_1,x_2,\dots,x_k]^q=1$ for all 
$x_1,x_2,\dots,x_k\in G$. Then $\gamma_k(G)$ is locally finite.
\bigskip

This extends the positive solution of the Restricted Burnside 
Problem for groups of prime-power exponent (take $k=1$). 
Another result considered in Section 3 deals with residually finite groups 
in which the commutators $[x_1,x_2,\dots,x_k]$ are Engel. 
\bigskip

\noindent{\bf Theorem 3.2.} Let $k,n$ be positive integers, and $G$ a residually
finite group such that $[x_1,x_2,\dots,x_k]$ is $n$-Engel for any
$x_1,x_2,\dots,x_k\in G$. Then $\gamma_k(G)$ is locally nilpotent.
\bigskip

For $k=1$ this is a well-known result of J. Wilson \cite{w}.

Section 4 is devoted to some questions on exponent of a finite group with
automorphisms. The main result of the section is the following theorem.
\bigskip

\noindent{\bf Theorem 4.3} (Khukhro, Shumyatsky \cite{khsh}). Let $q$ be a 
prime, $n$ an integer.
Suppose that a non-cyclic group $A$ of order $q^2$ acts on a finite group 
$G$ of coprime order in such a manner that the exponents of the 
centralizers $C_G(a)$ of non-trivial elements $a\in A$ divide $n$. 
Then the exponent of $G$ is bounded in terms of $n$ and $q$.
\bigskip

Note that the exponent of the centralizer of a single automorphism $a$
of a finite group $G$ has no impact over the exponent of $G$. Indeed, any abelian
group of odd order admits a fixed-point-free automorphism of order two. Hence,
we cannot bound the exponent of $G$ solely in terms of the exponent of
$C_G(a)$. In view of this the following theorem seems to be interesting.
\bigskip

\noindent{\bf Theorem 4.4.}
Let $n$ be a positive integer, $G$ a finite group of odd order admitting an involutory automorphism $a$ such that $C_G(a)$ is of exponent dividing $n$.
Suppose that for any $x\in G$ the element $[x,a]=x^{-1}x^a$
has order dividing $n$. Then the exponent of $G$ is bounded in terms of $n$.
\bigskip

In Section 5 we obtain some sufficient conditions for a periodic residually 
finite group to be locally finite. 
\bigskip

\noindent{\bf Theorem 5.3} (Shalev \cite{sha4}).
Let $G$ be a periodic residually finite group having a finite 2-subgroup $A$ 
such that the centralizer $C_G(A)$ is finite. Then $G$ is locally finite.  
\bigskip

\noindent{\bf Theorem 5.4.} 
Let $q$ be a prime, $G$ a residually finite group in which each 2-generator 
subgroup is finite. Suppose that $G$ has a finite $q$-subgroup $A$ 
such that the centralizer $C_G(A)$ is finite. Then $G$ is locally finite.  
\bigskip

This paper certainly is not a comprehensive survey on Lie methods in group 
theory. Some important areas are not even mentioned here. For the reader 
willing to learn more on this subject we recommend the survey \cite{sha2}
and the textbooks \cite[Chapter VIII]{hb}, \cite{khu1}, \cite{vl}.
Throughout the paper $p$ stands for a fixed prime. We use the term\linebreak
``$\{a,b,c,\dots\}$-bounded" to mean ``bounded from above by some function of
$a,b,c,\dots$".

I am grateful to A. Mann for a number of suggestions simplifying proofs in Section 3.
\bigskip

\bigskip

\bigskip

\noindent{\bf 2. Associating a Lie algebra to a group}
\bigskip

Let $L$ be a Lie algebra over a field ${\mathfrak k}$.
Let $k,n$ be positive integers and let $x_1,x_2,\dots,x_k,\newline x,y$ be elements 
of $L$. We define inductively
$$[x_1]=x_1;\ [x_1,x_2,\dots,x_k]=[[x_1,x_2,\dots,x_{k-1}],x_k]$$
and
$$[x,{}_0y]=x;\ [x,{}_ny]=[[x,{}_{n-1}y],y].$$
An element $a\in L$ is
called ad-nilpotent if there exists a positive integer $n$ such that
$[x,{}_na]=0$ for all $x\in L$. If $n$ is the least 
integer with the above property then we say that $a$ is ad-nilpotent of index 
$n$. Let $X\subseteq L$ be any subset of $L$. By a commutator in elements of $X$
we mean any element of $L$ that can be obtained as a Lie product of elements of
$X$ with some system of brackets. Denote by $F$ the free Lie algebra over 
${\mathfrak k}$ on countably many free generators $x_1,x_2,\dots$. Let
$f=f(x_1,x_2,\dots,x_n)$ be a non-zero element of $F$. The algebra $L$ is said
to satisfy the identity $f\equiv 0$ if $f(a_1,a_2,\dots,a_n)=0$ for any
$a_1,a_2,\dots,a_n\in L$. In this case we say that $L$ is PI. The next
theorem is straightforward from Theorem 1.5 quoted in the introduction.
\bigskip

\noindent {\bf Theorem 2.1.} Let $L$ be a Lie algebra over a field
${\mathfrak  k}$ generated by 
$a_1,a_2,\dots,a_m$. Assume that $L$ satisfies an identity $f\equiv 0$ and that
each commutator in the generators $a_1,a_2,\dots,a_m$ is ad-nilpotent of index 
at most $n$. Then $L$ is nilpotent of $\{f,n,m,{\mathfrak k}\}$-bounded class.
\medskip

\noindent{\bf Proof.}
Consider the free $m$-gene\-ra\-ted 
Lie ${\mathfrak k}$-algebra $F_m$ on free gene\-ra\-tors $f_1,\ldots ,f_m$, 
and let $T$ be the ideal of $F_m$ gene\-ra\-ted by all 
values of $f$ on the elements of $F_m$ and by all elements of the form 
$[g,{}_nc]$, where $g\in F_m$ and $c$ is an 
arbitrary commutator in the $f_i$. Then $F/T$ satisfies the 
hypothesis of Theorem 1.5 and hence is nilpotent of some class $u=u(m,n,f,\mathfrak k)$. 
Clearly, $L$ is the image of $F/T$ under the homomorphism 
induced by the mapping $f_i\rightarrow a_i$.
Hence $L$ is nilpotent of class at most $u$. $\Box$
\bigskip

An important criterion for a Lie algebra to be PI is the following
\bigskip

\noindent {\bf Theorem 2.2} (Bahturin-Linchenko-Zaicev).
Let $L$ be a Lie algebra over a field ${\mathfrak k}$. Assume that a finite 
group $A$ acts on $L$ by automorphisms in such a manner that $C_L(A)$, the
subalgebra formed by fixed elements, is PI. Assume further that 
the characteristic of ${\mathfrak k}$ is either 0 or prime to the order of $A$. 
Then $L$ is PI.
\bigskip

This theorem was first proved by Yu.\,A.\,Bahturin and 
M.\,V.\,Zaicev for solvable 
groups $A$ \cite{bz} and later extended by V.\,Linchenko to the general case \cite{li}.
\bigskip

\noindent {\bf Corollary 2.3.} 
Let $F$ the free Lie algebra of countable rank over ${\mathfrak k}$. Denote by 
$F^*$ the set of non-zero elements of $F$. For any finite group $A$ there 
exists a mapping $$\phi:F^*\rightarrow F^*$$ such that if $L$ and $A$ are
as in Theorem 2.2, and if $C_L(A)$ satisfies an identity $f\equiv 0$, then $L$ 
satisfies the identity $\phi(f)\equiv 0$.
\medskip

\noindent{\bf Proof.} Assume that the assertion is false. Then there exists
$f\in F^*$ such that for any $g\in F^*$ we can choose a Lie algebra $L_g$ which
does not satisfy the identity $g\equiv 0$ and admits an action by $A$ with
$C_{L_g}(A)$ satisfying $f\equiv 0$. Consider the direct sum
$M=\oplus L_g$, where the summation is taken over all $g\in F^*$. It is easy
to see that $M$ also admits an action by $A$ with $C_M(A)$ satisfying
$f\equiv 0$. Obviously $M$ cannot be PI, a contradiction against 2.2.
$\Box$
\bigskip

We now turn to groups. Let $G$ be any group. For $x,y\in G$ we
use $[x,y]$ to denote the group-commutator $x^{-1}y^{-1}xy$.
The long commutators $[x_1,x_2,\dots,x_k]$ and $[x,{}_ny]$ are 
defined as in Lie algebras:
$$[x_1]=x_1;\ [x_1,x_2,\dots,x_k]=[[x_1,x_2,\dots,x_{k-1}],x_k]$$
and
$$[x,{}_0y]=x;\ [x,{}_ny]=[[x,{}_{n-1}y],y].$$
If $x_1,x_2,\dots,x_k$ belong to a set $A$ then we say that
$[x_1,x_2,\dots,x_k]$ is a simple commutator of weight $k$ in elements of $A$. 

The following commutator identities hold in any group and can be checked 
manually.
\bigskip

\noindent{\bf Lemma 2.4.}   $$[x,y]^{-1}=[y,x]$$
   $$[xy,z]=[x,z][x,z,y][y,z]$$
   $$[x,yz]=[x,z][x,y][x,y,z]$$
   $$[x,y^{-1},z]^y[y,z^{-1},x]^z[z,x^{-1},y]^x=1.$$
\bigskip

Let us also record the Collection Formula (see \cite[p. 240]{hb}).
For any integer $n$ and any subgroup $H$ of a group $G$ we denote by $H^n$
the subgroup of $G$ generated by the $n$-th powers of elements from $H$
and use $\gamma_n(G)$ for the n-th term of the lower central series of $G$. 
\bigskip

\noindent{\bf Lemma 2.5.} Let $x,y$ be elements of a group $G$. Then
$$(xy)^{p^n}\equiv x^{p^n}y^{p^n}\ {mod}\
 \gamma_2(G)^{p^n}\prod\limits_{1\leq r\leq n}
\gamma_{p^r}(G)^{p^{n-r}}.$$
\bigskip

A series of subgroups $$G=G_1\geq G_2\geq\dots   \eqno{(*)}$$
is called an $N$-series if it satisfies $[G_i,G_j]\leq G_{i+j}$ for all $i,j$.
Obviously any $N$-series is central, i.e. $G_i/G_{i+1}\leq Z(G/G_{i+1})$ for any
$i$. An $N$-series is called $N_p$-series if $G_i^p\leq G_{pi}$ for all $i$.

Generalizing constructions discovered by Magnus \cite{ma} and Zassenhaus 
\cite{za}, Lazard noticed in 
\cite{la} that a Lie ring $L^*(G)$ can be associated to
any $N$-series $(*)$ of a group $G$. He also discovered some very useful
special properties of $L^*(G)$ in case $(*)$ is an $N_p$-series. Let us briefly 
describe the construction. 

Given an $N$-series $(*)$, let $L^*(G)$ be the direct sum of the abelian groups
$L_i^*=G_i/G_{i+1}$, written additively. Commutation in $G$ induces a binary
operation $[,]$ in $L$. For homogeneous elements 
$xG_{i+1}\in L_i^*, yG_{j+1}\in L_j^*$ 
the operation is defined by 
$$[xG_{i+1},yG_{j+1}]=[x,y]G_{i+j+1}\in L_{i+j}^*$$ and extended to arbitrary 
elements of $L^*(G)$ by linearity. It is easy to check that the operation is 
well-defined and that $L^*(G)$ with the operations $+$ and $[,]$ is a Lie ring. 

The above procedure can be performed for each $N$-series of $G$. If all
quotients $G_i/G_{i+1}$ of an $N$-series $(*)$ have exponent $p$
then $L^*(G)$ can be viewed as a Lie algebra
over $\mathbb F_p$, the field with $p$ elements.
This is always the case if $(*)$ is an $N_p$-series.
We are now concerned with the relationship between $G$ and $L^*(G)$. For any
$x\in G_i\setminus G_{i+1}$ let $x^*$ denote the element $xG_{i+1}$ of $L^*(G)$.
\bigskip

\noindent{\bf Lemma 2.6} (Lazard, \cite{la}). If $(*)$ is an $N_p$-series then $(ad\, x^*)^p=ad\, (x^p)^*$ for any $x\in G$. Consequently, 
if $x$ is of finite order $t$ then $x^*$ is ad-nilpotent of index at most $t$.

\bigskip

Let $Fr$ denote the free group on free generators $x_1,x_2,\dots$, and choose
a non-trivial element $w=w(x_1,x_2,\dots,x_s)\in Fr$. 
We say that a group $G$ satisfies the identity $w\equiv 1$ if 
$w(g_1,g_2,\dots,g_s)=1$ for any $g_1,g_2,\dots,g_s\in G$.

The following proposition can be deduced from the proof of Theorem 1 in the 
paper of J. Wilson and E. Zelmanov \cite{wize}
\bigskip

\noindent{\bf Proposition 2.7.} Let $G$ be a group satisfying an identity
$w\equiv 1$. Then there exists a non-zero Lie polynomial $f$ over $\mathbb F_p$ 
depending only on $p$ and $w$ such that for any $N_p$-series $(*)$ of $G$ the
algebra $L^*(G)$ satisfies the identity $f\equiv 0$.
\bigskip

In fact J. Wilson and E. Zelmanov describe in \cite{wize} an effective algorithm
allowing to write $f$ explicitely for any $p$ and $w$ but we do not require 
this. Let us just record an important special case of the above proposition.
\bigskip

\noindent{\bf Proposition 2.8} (Higman, \cite{hi}). Let $n$ be a $p$-power,
$G$ a group such that $x^n=1$ for any $x\in G$. 
Then for any $N_p$-series $(*)$ the algebra $L^*(G)$ satisfies the identity
$$\sum\limits_{\pi\in S_{n-1}}[x_0,x_{\pi(1)},x_{\pi(2)},\ldots,x_{\pi 
(n-1)}]=0.\eqno{(2.9)}$$ 
\bigskip

In general a group $G$ has many $N_p$-series so that there are many ways to
associate to $G$ a Lie algebra as described above. We will introduce now an
$N_p$-series which is of particular importance for applications of Lie-theoretic
results to group theory.

To simplify the notation we write $\gamma_i$ for $\gamma_i(G)$.
Set $D_i=D_i(G)=\prod\limits_{jp^k\geq i}\gamma_j^{p^k}$.
The subgroups $D_i$ form a series $G=D_1\geq D_2\geq\dots$ in the group $G$.
\bigskip

\noindent{\bf Proposition 2.10} (\cite[p. 250]{hb}). The series $\{D_i\}$ is an
$N_p$-series.
\bigskip

We will call $\{D_i\}$  the $p$-dimension central series of $G$. It is also 
known as the Lazard series or the Jennings-Lazard-Zassenhaus series.

In conformity with the earlier described procedure we can 
associate to $G$ a Lie algebra $DL(G)=\oplus L_i$ over $\mathbb F_p$ corresponding 
to the $p$-dimension central series. Here $L_i=D_i/D_{i+1}$. This algebra plays a 
crucial r\^ole in all results considered in the paper.

Let $L_p(G)=\langle L_1\rangle$ be the subalgebra of $DL(G)$ generated by $L_1$.
If $G$ is finitely generated then nilpotency of $L_p(G)$ has strong impact over
the structure of $G$.  The following proposition is implicit in \cite{ze3}.
\bigskip

\noindent{\bf Proposition 2.11.} Let $G$ be generated by elements 
$a_1,a_2,\dots,a_m$, 
and assume that $L_p(G)$ is nilpotent of class at most $c$.
Let $\rho_1,\rho_2,\dots,\rho_s$ be the list of all simple commutators in
$a_1,a_2,\dots,a_m$ of weight $\leq c$. Then for any non-negative integer $i$ 
the group $G$ can be written as a product 
$$G=\langle\rho_1\rangle\langle\rho_2\rangle\dots\langle\rho_s\rangle D_{i+1}$$
of the cyclic subgroups generated by
$\rho_1,\rho_2,\dots,\rho_s$ and $D_{i+1}$.  
\medskip

\noindent{\bf Proof.} We start with the following remark. For any positive 
integer $i$ the subgroup $D_i$ is generated by $D_{i+1}$ and elements of the
form 
$[b_1,\dots,b_j]^{p^k}$,
where $jp^k\geq i$ and $b_1,\dots,b_j\in\{a_1,\dots,a_m\}$.
This can be shown using for example formulae 2.4 and 2.5.

The proposition will be proved by induction on $i$, the case $i=0$ being 
trivial. Assume that $i\geq 1$ and
$$G=\langle\rho_1\rangle\langle\rho_2\rangle\dots\langle\rho_s\rangle D_{i}.$$
Then any element $x\in G$ can be written in the form
$$x=\rho^{\alpha_1}_1\rho^{\alpha_2}_2\dots\rho^{\alpha_s}_sy,\eqno (2.12)$$
where $y\in D_i$. Without any loss of generality we can assume that $D_{i+1}=1$.

By the remark made in the beginning of the proof we can write
$$y=
{(\sigma_1^{p^{k_1}})}^{\beta_1} 
{(\sigma_2^{p^{k_2}})}^{\beta_2}\dots 
{(\sigma_t^{p^{k_t}})}^{\beta_t},\eqno (2.13)$$ where each 
$\sigma_n$
is of the form 
$[b_1,\dots,b_j]$,
with $jp^{k_n}\geq i$ and 
$b_1,\dots,b_j\in\{a_1,\dots,a_m\}$.

Let ${\tilde a_l}$ denote $a_lD_2\in L_p(G);\ l=1,\dots,m$.
By the hypothesis $L_p(G)$ is nilpotent of class $c$, that is
$[{\tilde b_1},\dots,{\tilde b_{c+1}}]=0$ for any
$b_1,\dots,b_{c+1}\in\{a_1,\dots,a_m\}$. This implies that 
$[b_1,\dots,b_{c+1}]\in D_{c+2}$ for any
$b_1,\dots,b_{c+1}\in\{a_1,\dots,a_m\}$ and 
$\gamma_{c+1}\leq D_{c+2}$. 
Then, by Proposition 2.10, for any $d\geq c+1$ we have 
$\gamma_{d}\leq D_{d+1}$. 

Now, if $\sigma_n$ is of the form $[b_1,\dots,b_j]$ with $j\geq c+1$ then
$$\sigma_n^{p^{k_n}}\in\gamma_j^{p^{k_n}}\leq D_{j+1}^{p^{k_n}}\leq 
D_{(j+1)p^{k_n}}\leq D_{i+1}=1.$$ Hence we can assume that each 
$\sigma_n$ is of the form $[b_1,\dots,b_j]$ with $j\leq c$, in which case
$\sigma_n$ belongs to the list $\rho_1,\rho_2,\dots,\rho_s$.

It remains to remark that by 2.10 $\sigma_n^{p^{k_n}}\in Z(G)$. Comparing now
(2.12) and (2.13) we obtain that
$$x\in\langle\rho_1\rangle\langle\rho_2\rangle\dots\langle\rho_s\rangle,$$
as required.
$\Box$

\bigskip
The following corollary is now immediate.
\bigskip

\noindent{\bf Corollary 2.14.} Assume the hypotheses of 2.6 and suppose that each
$\rho_j$ is of order at most $K$. Then $D_i$ is of index at most $K^s$ for any 
$i$. In particular, the series $D_i$ becomes stationary after finitely many 
steps. 
\bigskip

A group $G$ is said to be residually-p if for any non-trivial element $x\in G$ 
there exists a normal subgroup $N\leq G$ such that $x\not\in N$ and $G/N$ is a 
finite $p$-group.  M. Lazard proved in \cite{l2} that if $G$ is a finitely 
generated pro-p group such that $L_p(G)$ is nilpotent then $G$ is $p$-adic 
analytic. Since any finitely generated residually-p group can be embedded in a 
finitely generated pro-p group (see for example \cite{ddms}), we obtain the 
following 
\bigskip

\noindent{\bf Proposition 2.15.} 
If $G$ is a finitely generated residually-p group
such that $L_p(G)$ is nilpotent then $G$ has a faithful linear representation
over the field of $p$-adic numbers.
\bigskip

\bigskip

\bigskip

\noindent{\bf 3. Groups with commutators of bounded order}
\bigskip

This section is devoted to groups in which all commutators
$[x_1,x_2,\dots,x_k]$ satisfy some restrictive condition.
In particular we describe a proof of the following theorem. 
\bigskip

\noindent{\bf Theorem 3.1.}
Let $k$ be an integer, $q=p^s$ a prime-power, $G$ a residually finite group such
that $[x_1,x_2,\dots,x_k]^q=1$ for all 
$x_1,x_2,\dots,x_k\in G$. Then $\gamma_k(G)$ is locally finite.
\bigskip

This generalizes the positive solution of the Restricted Burnside 
Problem for groups of prime-power exponent (take $k=1$). 
The theorem is no longer true if the assumption that $G$ is residually finite
is dropped.
Using the technique developed by A. Ol'shanskii \cite{ol} it is possible (for 
sufficiently big values of $q$) to construct a group $G$ that satisfies the 
identity $[x,y]^q=1$ and has the derived group $G'$
non-periodic.

Another result considered in this section deals with residually finite groups 
in which the commutators $[x_1,x_2,\dots,x_k]$ are Engel. An element $x$ of a 
group $G$ is called  (left) $n$-Engel if $[g,{}_nx]=1$ for any
$g\in G$. A group $G$ is called $n$-Engel if all elements of $G$ are $n$-Engel. 
It is a long-standing problem whether any $n$-Engel group is locally nilpotent.
In \cite{w}  J. Wilson proved that this is true if $G$ is residually finite.
His proof relies on a result of A. Shalev \cite{sha1} which uses the positive
solution of the Restricted Burnside Problem \cite{ze1},\cite{ze2}. We will prove
\bigskip

\noindent{\bf Theorem 3.2.} Let $k,n$ be positive integers, $G$ a residually
finite group such that $[x_1,x_2,\dots,x_k]$ is $n$-Engel for any
$x_1,x_2,\dots,x_k\in G$. Then $\gamma_k(G)$ is locally nilpotent.
\bigskip

We will consider first finite groups satisfying the hypotheses of 3.1 and 3.2. 
To prepare the use of Theorem 1.5 we will show that if $G$ is a finite group
satisfying the hypothesis of Theorem 3.1 (respectively, 3.2),
then $\gamma_k(G)$ is a $p$-group (respectively, is nilpotent).
\bigskip

In the proof of the next lemma we follow advice of A. Mann and R. Solomon.
\bigskip

\noindent {\bf Lemma 3.3.} Let $G$ be a finite group in which all commutators
$[x_1,x_2,\dots,x_k]$ are $p$-elements. Then $\gamma_k(G)$ is a $p$-group.
\medskip

\noindent {\bf Proof.} Assume that the result is false and let $G$ be a 
counterexample of minimal possible order. Obviously $G$ has no non-trivial
normal $p$-subgroups.
Let $r$ be a prime divisor of $|G|$ distinct from $p$ and suppose that $G$ contains an $r$-subgroup $R$ such that
the quotient $N_G(R)/C_G(R)$ has an element $z$ of order prime to $r$. Then $z$ can be viewed as a non-trivial automorphism
of $R$ and therefore
$[R,\underbrace{z,\dots,z}_{k-1}]\neq 1$, \cite[Theorem 5.3.6]{go}. On the other hand, by the hypothesis, 
$[R,\underbrace{z,\dots,z}_{k-1}]$ must contain $p$-elements,
a contradiction. Therefore $N_G(R)/C_G(R)$ is an $r$-group for
any $r$-subgroup $R$ of $G$. Hence $G$ possesses a normal $r$-complement $K$ \cite[7.4.5]{go}. 
The induction on $|G|$ shows
that $\gamma_k(K)$ is a $p$-group.
Since $G$ has no non-trivial normal $p$-subgroups, it follows that $\gamma_k(K)=1$ and $K$ is nilpotent. Again because 
$G$ has no non-trivial normal $p$-subgroups, we conclude that
$K$ is a $p'$-group. But then so is $G$, a contradiction.
$\Box$ 
\bigskip

\noindent{\bf Lemma 3.4.} If for any elements 
$y,x_1,x_2,\dots,x_k$ of a finite group $G$ there exists an integer $n$ such 
that $[y,{}_n[x_1,x_2,\dots,x_k]]=1$ then 
$\gamma_k(G)$ is nilpotent.
\medskip

\noindent {\bf Proof.} Engel elements of any finite group generate a nilpotent subgroup \cite[III, 6.14]{hu}.
Since $G$ is finite and $\gamma_k(G)=
\langle [x_1,x_2,\dots,x_k];\ x_1,x_2,\dots,x_k\in G\rangle$, the result 
follows. 
$\Box$ 
\bigskip

The next lemma is a particular case of Lemma 2.1 in the paper of 
J. Wilson \cite{w}.
\bigskip

\noindent{\bf Lemma 3.5.} Let $G$ be a finitely generated residually
finite-nilpotent group. For each prime $p$ let $J_p$ be the intersection of all
normal subgroups of $G$ of finite $p$-power index. If $G/J_p$ is nilpotent
for each $p$ then $G$ is nilpotent.
\medskip

\noindent{\bf Proof.} Let $R$ be the intersection of all normal subgroups $N$ of
$G$ such that $G/N$ is torsion-free and nilpotent. Since a torsion-free 
nilpotent group is residually-$p$ for any prime $p$, it follows that $J_p\leq R$
for any $p$. We therefore conclude that $G/R$ is nilpotent. Then so is
$G/[R,\underbrace{G,\dots,G}_{k}]$ for any positive integer $k$. Set $R_0=R$,
$R_{k+1}=[R,\underbrace{G,\dots,G}_{k}]$. Since $G$ is residually nilpotent,
the intersection $\cap_iR_i$ is trivial. By the choice of $R$ the quotient
$R/R_1$ must be periodic. This is a subgroup of a finitely generated nilpotent
group $G/R_1$ and therefore $R/R_1$ is finitely generated too. Hence $R/R_1$
is finite, and let $\pi$ be the set of primes dividing the order of $R/R_1$. 
Then $R/R_k$ is a $\pi$-group for any $k$. (It is a well-known general fact:
If $R$ is any normal subgroup of a nilpotent group $G$, and if $R/[R,G]$ is a 
$\pi$-group, then $R$ is a $\pi$-group.) It follows that $\cap_{p\in\pi}J_p\leq
\cap_iR_i=1$. Since $\pi$ is finite, $G$ embeds into the direct product of
finitely many nilpotent groups $G/J_p,\ p\in\pi$. Hence $G$ is nilpotent. $\Box$
\bigskip

\noindent {\bf Proof of Theorem 3.1.} 
By Lemma 3.3 $\gamma_k(G)$ is residually-p. 
To prove that $\gamma_k(G)$ is locally finite let us take arbitrarily a 
finite subset $S$ of $\gamma_k(G)$ and show that $S$ generates a finite 
subgroup. Since $\gamma_k(G)$ is generated by elements of the form
$[x_1,x_2,\dots,x_k]$, it follows 
that $S$ is contained in some subgroup $H$ generated by finitely many elements 
$a_1,a_2,\dots,a_m$ each of which can be written in the form
$[x_1,x_2,\dots,x_k]$. Now it is important to observe that if $\rho$ is an 
arbitrary commutator in $a_1,a_2,\dots,a_m$ with some system of brackets then
$\rho$ also can be written in the form
$\rho=[x_1,x_2,\dots,x_k]$ for suitably chosen $x_1,x_2,\dots,x_k\in G$. 
We want to show that $H$ is finite.
Let $L=L_p(H)$ be the Lie algebra associated with the $p$-dimension central series 
\begin{center}
$H=D_1\geq D_2\geq\dots$
\end{center}
of $H$ (see Section 2).  Then $L$ is generated by
$\tilde a_i=a_iD_2;\ i=1,2,\dots m$. Let $b$ be any Lie-commutator in
$\tilde a_1, \tilde a_2,\dots, \tilde a_m$ and $c$ the group-commutator in
$a_1,a_2,\dots,a_m$ having the same system of brackets as $b$. We have already
observed that $c=[x_1,x_2,\dots,x_k]$ for some $x_1,x_2,\dots,x_k\in G$, so we
conclude that $c^q=1$. Suppose that $c\in D_j\setminus D_{j+1}$. It is immediate from
the definition of $L$ that either $b=0$ or $b=cD_{j+1}$. This implies 
(Lemma 2.6) that $b$ is ad-nilpotent of index at most $q$.
Further, the group $H$ satisfies the identity $[y_1,y_2,\dots,y_k]^q=1$.
Therefore, by 2.7, $L$ satisfies some non-trivial polynomial identity. 
By Theorem 1.5 we conclude that $L$ is nilpotent. Corollary 2.14 now shows
that the index of $D_i$ in $H$ is bounded from above by some number that 
depends only on $q$, $m$ and the nilpotency class of $L$. Since $H$ is 
residually-p, it follows that the intersection of all $D_i$ is trivial
 and we derive that $H$ is finite.
$\Box$
\bigskip

The proof of 3.2 is similar to that of 3.1.
\bigskip

\noindent {\bf Proof of Theorem 3.2.} 
By Lemma 3.4 $\gamma_k(G)$ is residually finite-nilpotent. 
To prove that $\gamma_k(G)$ is locally nilpotent let us take arbitrarily a 
finite subset $S$ of $\gamma_k(G)$ and show that $S$ generates a nilpotent 
subgroup. Since $\gamma_k(G)$ is generated by elements of the form
$[x_1,x_2,\dots,x_k]$, it follows 
that $S$ is contained in some subgroup $H$ generated by finitely many elements 
$a_1,a_2,\dots,a_m$ each of which can be written in the form
$[x_1,x_2,\dots,x_k]$. As in the proof of 3.1 we observe that if $\rho$ is an 
arbitrary commutator in $a_1,a_2,\dots,a_m$ with some system of brackets then
$\rho$ also can be written in the form
$\rho=[x_1,x_2,\dots,x_k]$ for suitably chosen $x_1,x_2,\dots,x_k\in G$. 
We want to show that $H$ is nilpotent.
Lemma 3.5 allows us to assume that $H$ is residually-p for some prime $p$. 
Let $L=L_p(H)$ be the Lie algebra associated with the $p$-dimension central series 
\begin{center}
$H=D_1\geq D_2\geq\dots$
\end{center}
of $H$.  Then $L$ is generated by
$\tilde a_i=a_iD_2;\ i=1,2,\dots m$. Let $b$ be any Lie-commutator in
$\tilde a_1, \tilde a_2,\dots, \tilde a_m$ and $c$ be the group-commutator in
$a_1,a_2,\dots,a_m$ having the same system of brackets as $b$. We have already
observed that $c=[x_1,x_2,\dots,x_k]$ for some $x_1,x_2,\dots,x_k\in G$, so we
conclude that $c$ is $n$-Engel. This implies that $b$ is ad-nilpotent of index
at most $n$.
Further, the group $H$ satisfies the identity $[y,{}_n[y_1,y_2,\dots,y_k]]=1$.
Therefore, by 2.7, $L$ satisfies some non-trivial polynomial identity. 
Theorem 1.5 now implies that $L$ is nilpotent. Hence, by
Proposition 2.15, $H$ has a faithful linear representation over the field of 
$p$-adic numbers. Clearly $H$ cannot have a free subgroup of rank two, and so 
by Tits' Alternative \cite{tits} $H$ has a solvable
subgroup of finite index. Lemma 3.4 shows that each finite quotient of $H$ is 
solvable, so that we conclude that $H$ is solvable too. It is now easy to 
deduce from a result of K. Gruenberg \cite{gr} that $H$ is nilpotent, 
as required.
$\Box$ 
\bigskip

We conclude this section by listing some open problems related to the results
described here.
\bigskip

\noindent{\bf Problem 1.} Does there exist a group $G$ satisfying the 
hypothesis of Theorem 3.1 and having $\gamma_k(G)$ of unbounded exponent?
\bigskip

The answer to this problem is likely to be ``yes" but finding an example does 
not seem to be easy. Using results of B. Hartley \cite{ha4} one
can show that such an example cannot be finitely generated (
see \cite{shu0} for detail in the case $k=2$).
\bigskip

\noindent{\bf Problem 2.}
Let $k,n$ be positive integers, $G$ a residually finite group such that\linebreak 
$[x_1,x_2,\dots,x_k]^n=1$ for any $x_1,x_2,\dots,x_k\in G$. Is then 
$\gamma_k(G)$ necessarily locally finite?
\bigskip

\noindent{\bf 4. Bounding the exponent of a finite group  with automorphisms} 
\bigskip

Let $q$ be a prime, and let $A$ be a non-cyclic group of order $q^2$ acting
on a finite group $G$. It is well-known (see \cite[Theorem 6.2.4, Theorem 5.3.16
]{go}) that if $G$ is any group of order prime to $q$ then
$$G=\left<C_G(a); a\in A^{\#}\right>,\eqno{(4.1)}$$ 
where $A^{\#}=A\setminus \{ 1 \}$. 
If, moreover, $G$ is a $p$-group then we even have 
$$G=\prod\limits_{a\in A^{\#}}C_G(a).\eqno{(4.2)}$$ 
Since $A$ normalizes some Sylow $p$-subgroup of $G$ for any $p$ dividing $|G|$
(\cite[Theorem 6.2.4]{go}), it is immediate that if $|C_G(a)|\leq n$ for any 
$a\in A^\#$ then 
$|G|\leq n^{q+1}$ (we use that $A$ has exactly $q+1$ cyclic subgroups). How 
profound is the connection between the structure of $G$ and that of $C_G(a);\
a\in A^\#$? It is known that if $C_G(a)$ is nilpotent for each $a\in A^\#$
then $G$ is metanilpotent \cite{wa}. In some
situations this result holds even if $G$ is allowed to be infinite periodic
\cite{shu0}.
In this section we will prove the following theorem.
\bigskip

\noindent{\bf Theorem 4.3} (Khukhro, Shumyatsky \cite{khsh}).
Suppose that $A$ is a non-cyclic group of order $q^2$ acting on a finite group 
$G$ of coprime order, and let $n$ be such an integer that the exponents of the 
centralizers $C_G(a)$ of non-trivial elements $a\in A^{\#}$ divide $n$. 
Then the exponent of $G$ is bounded in terms of $n$ and $q$.
\bigskip

Note that the exponent of the centralizer of a single automorphism $a$
of a finite
group $G$ has no impact over the exponent of $G$. Indeed, any abelian
group of odd order admits a fixed-point-free automorphism of order two. Hence,
we cannot bound the exponent of $G$ solely in terms of the exponent of
$C_G(a)$. In view of this
the following theorem seems to be interesting.
\bigskip

\noindent{\bf Theorem 4.4.}
Let $n$ be a positive integer, $G$ a finite group of odd order admitting an involutory automorphism $a$ such that $C_G(a)$ is of exponent dividing $n$.
Suppose that for any $x\in G$ the element $[x,a]=x^{-1}x^a$
has order dividing $n$.
Then the exponent of $G$ is bounded in terms of $n$.
\bigskip

Apart from the technique described in Section 2 the proof of Theorems 
4.3 and 4.4 involves using powerful $p$-groups. These were introduced by A.
Lubotzky and A. Mann in \cite{lbmn}: a finite $p$-group $G$ is powerful if and 
only if $G^p\geq [G,G]$ for $p\ne 2$ (or $G^4\geq [G,G]$ for $p=2$). 

Powerful $p$-groups have many nice linear properties, of which we need 
the following: if a powerful $p$-group $G$ is generated by elements of 
exponent $p^e$, then the exponent of $G$ is $p^e$ too (see 
\cite[Lemma 2.2.5]{ddms}). Combining this with 4.1, we conclude
that if $G$ is a powerful $p$-group satisfying the hypothesis of Theorem 4.3, 
then the exponent of $G$ divides $n$. Thus, it is sufficient to reduce the 
proof of Theorem 4.3 to the case of powerful $p$-groups. A. Shalev was the 
first to discover relevance of powerful $p$-groups to problems on automorphisms 
of finite groups \cite{sha5}, \cite{sha6}. However our situation is quite 
different from that considered in the papers of A. Shalev. The reason powerful 
$p$-groups emerge in the context of the Restricted Burnside Problem is the following lemma.
\bigskip

\noindent{\bf Lemma 4.5.} Suppose that $P$ is a $d$-generator finite 
$p$-group such that the Lie algebra $L_p(P)$ is nilpotent of class $c$. Then 
$P$ has a powerful characteristic subgroup of $\{p,c,d\}$-bounded index.
\medskip

\noindent{\bf Proof.} Let $\rho_1,\dots,\rho_s$ be all simple commutators of 
weight $\leq c$ in the generators of $P$; here $s$ is a $\{d,c\}$-bounded number. 
Since $P$ is a finite $p$-group, the nilpotency of $L_p(P)$ of class $c$ 
implies that every element $g\in P$ can be written in the form 
$g=\rho_1^{k_1}\dots\rho_s^{k_s}$ (see Proposition 2.11).
Hence $ |P/P^{p^m}|\leq p^{sm}$ for any positive integer $m$. 
Let $V$ be the intersection of the kernels of all 
homomorphisms of $P$ into $GL_s(\mathbb F_p)$. Put $W=V$ if $p\ne 2$ 
(or $W=V^2$ if $p=2$). The exponent of the Sylow $p$-subgroup of $GL_s(\mathbb
F_p)$ is a $\{p,s\}$-bounded number. Then $P^{p^a}\leq W$ for some 
$\{p,s\}$-bounded number $a$, which is also $\{p,c,d\}$-bounded, since $s$ is 
$\{c,d\}$-bounded. There is a $\{p,c,d\}$-bounded number $u\geq a$ such that 
$|P^{p^u}/P^{p^{u+1}}|\leq p^s$, for otherwise the inequality 
$|P/P^{p^m}|\leq p^{sm}$ would be violated for some $m$. Then 
$P^{p^u}\leq P^{p^a}\leq W$, and $P^{p^u}$ is $s$-generator since 
$|P^{p^u}/\Phi (P^{p^u})|\leq |P^{p^u}/P^{p^{u+1}}|\leq p^s$. Now, by 
\cite[Proposition 2.12]{ddms}, $P^{p^u}$ is a powerful subgroup. 
The index of $P^{p^u}$ is at most $p^{us}$ and hence is $\{p,c,d\}$-bounded.
$\Box$\bigskip

We will now quote a well-known lemma which is of fundamental importance whenever
one uses Lie algebra methods to study finite groups with automorphisms of
coprime order.
\bigskip

\noindent{\bf Lemma 4.6} (\cite[6.2.2 (iv)]{go}). Let $A$ be a finite group 
acting on a finite group $G$.  Assume that $(|A|,|G|)=1$ and that $N$ is a 
normal $A$-invariant subgroup of $G$. Then $C_{G/N}(A)=C_G(A)N/N$.
\bigskip

Let $H$ be a subgroup of a group $G$. Set $H_j=D_j\cap H$, where $D_j$ is the
$j$-th term of the $p$-dimension central series of $G$. Write
$$L(G,H)=\oplus H_jD_{j+1}/D_{j+1} \mbox{ and } L_p(G,H)=L_p(G)\cap L(G,H).$$

\noindent{\bf Observation 4.7.} $L(G,H)$ is a subalgebra of $DL(G)$. It is 
isomorphic to the Lie algebra associated with the $N_p$-series $\{H_j\}$ of $H$.
\bigskip

Let now a finite group $A$ act on a group $G$. Obviously $A$ induces
an automorphism group of every quotient $D_j/D_{j+1}$. This action extends to
the direct sum $\oplus D_j/D_{j+1}$. Thus, $A$ can be viewed as a group acting 
on $L_p(G)$ by Lie algebra automorphisms. The following remark is immediate 
from Lemma 4.6.
\bigskip

\noindent{\bf Observation 4.8.} If $G$ is finite and $(|G|,|A|)=1$ then 
$$L_p(G,C_G(A))=C_{L_p(G)}(A).$$
\bigskip

The next lemma will be helpful in the proof of Theorem 4.3.
\bigskip

\noindent{\bf Lemma 4.9.} Suppose that $L$ is a Lie algebra, $H$ 
a subalgebra of $L$ gene\-ra\-ted by $r$ elements $h_1,\dots ,h_r$ 
such that all commutators in the $h_i$ are ad-nilpotent in $L$ of index $t$. 
If $H$ is nilpotent of class $u$, then for some $\{r,t,u\}$-bound\-ed number $v$ 
we have $[L,\underbrace{H,\dots,H}_{v}]=0$. 
\medskip

\noindent{\bf Proof.} We apply to a 
sufficiently long (but of $\{r,t,u\}$-bound\-ed length) commutator $$[l, 
\,h_{i_1},\,h_{i_2},\,\ldots ]$$ a collecting process whose aim is to rearrange 
the $h_{i}$ (and emerging by the Jacobi identity commutators in the~$h_{i}$) 
in an ordered string after $l$, where all occurrences of a given element 
($h_i$ or a commutator in the $h_i$) would form an unbroken segment. Since $H$ 
is nilpotent, this process terminates at a linear combination of commutators 
with sufficiently long segments of equal elements. All these commutators 
are equal to $0$ because all commutators in the $h_i$ are ad-nilpotent by the 
hypothesis. $\Box$
\bigskip

\noindent{\bf Proof of Theorem 4.3.} 
As noted above, for every prime $p$ dividing $|G|$ there is an $A$-invar\-i\-ant
Sylow $p$-sub\-group~$P$ of~$G$. Since $P=\left< C_{P}(a); a\in A^{\#}\right>
$, such a prime $p$ must be a divisor of~$n$. If Theorem 4.3 is valid in the 
case where $G$ is a finite $p$-group, the exponents of all Sylow subgroups of 
$G$ are bounded in terms of $q$ and $n$, which implies that the exponent of $G$ 
is bounded. Thus, we may assume $G$ to be a finite $p$-group (for a prime 
$p\ne q$). We may also assume $G$ to be gene\-ra\-ted by $q^2$ elements, since 
every element $g\in G$ is contained in the $A$-inva\-ri\-ant subgroup~$\left< 
g^a; a\in A\right> $.
 
From now on, in addition to the hypothesis of Theorem 4.3, $G$ is a finite 
$q^2$-gener\-at\-or $p$-group and $n$ is a power of $p$. Let $A_1,A_2,\dots,
A_{q+1}$ be the distinct cyclic subgroups of $A$. Set $D_j=D_j(G)$, $L=L_p(G)$, $L_j=L\cap
D_j/D_{j+1}$, so that $L=\oplus L_j$. Let $L_{ij}=C_{L_j}(A_i)$. Then, by 4.2,
for any $j$ we have $$L_j=\sum\limits_{1\leq i\leq q+1}L_{ij}.$$
By Lemma 4.6 for any $l\in L_{ij}$ there exists $x\in D_j\cap C_G(A_i)$ such
that $l=xD_{j+1}$. By the hypothesis $x$ is of order at most $n$, whence $l$ is
ad-nilpotent of index at most $n$ (Lemma 2.6). Thus, $$\mbox{any element in }
L_{ij} \mbox{ is ad-nilpotent of index at most } n. \eqno (4.10)$$
Since $G$ is generated by $q^2$ 
elements, the $\mathbb F_p$-space $L_1$ is spanned by $q^2$ elements.
In particular, $L$ is gene\-ra\-ted by at most $q^2$ ad-nilpotent elements, each
from $L_{i1}$ for some $i$. But we cannot claim that every Lie commutator 
in these gene\-ra\-tors is again in some $L_{ij}$ and hence is ad-nilpotent too.

To overcome this difficulty, we extend the ground field of $L$ by 
a primitive $q$th root of unity~$\omega$, forming 
$\overline L=L\otimes \mathbb F_p [\omega]$. The idea is to replace $L$ by $\overline L$ 
and to prove that $\overline L$ is nilpotent of $\{q,n\}$-bound\-ed class, which will, of 
course, imply the same nilpotency result for $L$. Before that we translate 
the properties of $L$ into the language of $\overline L$. 

It is natural to identify $L$ with the $\mathbb F_p$-sub\-algebra 
$L\otimes 1$ of $\overline L$. We note that if an element $x\in L$ is ad-nilpotent 
of index $m$, say, then the ``same" element $x\otimes 1$ is ad-nilpotent in
$\overline L$ of the same index~$m$. 

Put 
$\overline{L_j}=L_j\otimes\mathbb F_p[\omega]$; then $\overline L=\left<\overline{L_1}\right>
$, since $L=\left<L_1\right>$, and $\overline L$ is the direct sum of the 
homogeneous components $\overline{L_j}$. 
Since the $\mathbb F_p$-space $L_1$ is spanned by $q^2$ elements, so is the $\mathbb F_p[\omega ]$-space $\overline{L_1}$.

The group $A$ acts naturally on $\overline L$, 
and  we have $\overline{L_{ij}}=C_{\overline{L_j}}(A_i)$, where
$\overline{L_{ij}}=L_{ij}\otimes\mathbb F_p[\omega]$.
Let us show that $$\mbox{any element } y\in\overline{L_{ij}} \mbox{ is 
ad-nilpotent of $\{q,n\}$-bounded index.}\eqno (4.11)$$
Since $\overline{L_{ij}}=L_{ij}\otimes\mathbb F_p[\omega]$, we can write
$$y=x_0+\omega x_1+\omega ^2x_2+ \dots +\omega ^{q-2}x_{q-2}$$ 
for some $x_s\in L_{ij}$, so that each of the summands $\omega^sx_s$ 
is ad-nilpotent of index $n$ by 4.10. Set $H=\left< 
x_0, \,\omega x_1, \,\dots ,\,\omega ^{q-2}x_{q-2}\right>$. Note 
that $H\subseteq C_{\overline L}(A_i)$, since $\omega^sx_s\in C_{\overline L}
(A_i)$ 
for all $s$. A commutator of weight $k$ in the $\omega^sx_s$ has the form
$\omega ^tx$ for some $x\in L_{im}$, where $m=kj$. By 4.10 such an $x$ is ad-nilpotent of 
index $n$ and hence so is $\omega ^tx$. 

Further, combining Observations 4.7 and 4.8 with Proposition 2.8, we conclude
that $C_L(A_i)$ satisfies the polynomial identity 2.9. This identity is
polylinear and so it is also satisfied by 
$C_L(A_i)\otimes\mathbb F_p[\omega]=C_{\overline L}(A_i)$. 
We have already observed that $H\subseteq C_{\overline L}(A_i)$, whence
the identity 2.9 is satisfied in $H$. Hence by Theorem 2.1 $H$ is
nilpotent of $\{q,n\}$-bound\-ed class. 
Lemma 4.9 now says that $[L,\underbrace{H,\dots,H}_{v}]=0$ for some 
$\{q,n\}$-bound\-ed number $v$. This establishes 4.11. 

Since $A$ is abelian, and the ground field is now a splitting field 
for $A$, every $\overline{L_j}$ decomposes in the direct sum of common 
eigenspaces for $A$. In particular, $\overline{L_1}$ is spanned by common 
eigenvectors for $A$, and it requires at most $q^2$ of them to span 
$\overline{L_1}$. Hence $\overline L$ is generated by $ q^2$ common 
eigenvectors for $A$ from $\overline{L_1}$. 
Every common eigenspace is contained in the centralizer 
$C_{\overline L}(A_i)$ for some $1\leq i\leq q+1$, since $A$ is 
non-cyclic. 
Note that any commutator in common eigenvectors is again a common 
eigenvector.
The main advantage of extending the ground field now becomes clear: 
if $l_1,\dots,l_{q^2}\in\overline{ L_1}$ are common eigenvectors for $A$
generating $\overline L$ then any commutator in these generators belongs to 
some $\overline{L_{ij}}$ and therefore, by 4.11, is ad-nilpotent of 
$\{q,n\}$-bounded index.

We already know that the identity 2.9 is satisfied in 
$C_L(A_i)\otimes\mathbb F_p[\omega]=C_{\overline L}(A_i)$. 
So, if $f$ denotes the Lie polynomial in 2.9 then, by 2.3, $\overline L$
satisfies some identity $\phi(f)\equiv 0$ which depends only on $n$ and $q$.
Theorem 2.1 now shows that $\overline L$ (hence $L$) is nilpotent of 
$\{q,n\}$-bounded class. 

By Lemma 4.5 $G$ contains a characteristic powerful subgroup $G_1$ of 
$\{q,n\}$-bounded index. Combining the hypothesis with 4.1 we see that $G_1$
is generated by elements of order dividing $n$. It follows that the exponent of 
$G_1$ divides $n$ \cite[2.2.5]{ddms}. Therefore the exponent of $G$ is 
$\{q,n\}$-bounded, as required. 
$\Box$
\bigskip

In the proof of Theorem 4.4 we use the following well-known fact.
\bigskip

\noindent{\bf Lemma 4.12.} Let $G$ be a finite group of odd order with an
involutory automorphism $a$. Then any element $x\in G$ can be written
uniquely in the form $x=gh$, where $g^a=g^{-1}$ and $h\in C_G(a)$. Moreover
$x^a=x^{-1}$ if and only if $x=[y,a]=y^{-1}y^a$ for some $y\in G$.
\bigskip

\noindent{\bf Proof of Theorem 4.4.} As in Theorem 4.3 we remark that $G$
possesses an $a$-invariant Sylow $p$-subgroup for any prime $p$ dividing $|G|$.
It is therefore sufficient to bound the exponent of $a$-invariant $p$-subgroups 
of $G$. So, without any loss of generality we may assume that $G$ is a $p$-group
and $n$ is a $p$-power. Take an arbitrary element $x\in G$. According to 4.12
we can write $x=gh$, where $g^a=g^{-1}$ and $h\in C_G(a)$. So to prove that
the order of $x$ is $n$-bounded it is sufficient to prove that the order
of $\langle g,h\rangle$ is $n$-bounded. Thus, we can assume that
$G=\langle g,h\rangle$.

Let $\{D_j\}$ be the $p$-dimension central 
series of $G$, $L_j=D_j/D_{j+1}$, $L=DL(G)$
the Lie algebra corresponding to the series $\{D_j\}$. We can naturally view
$a$ as an automorphism of $L$. Set $$L^+=C_L(a) \mbox{  and } L^-=\{l\in L;l^a=
-l\}.$$ Then one has:
$$[L^+,L^+]\leq L^+,\ [L^+,L^-]\leq L^-,\ [L^-,L^-]\leq L^+.\eqno (4.13)$$
For a fixed $j$ let us denote for a moment $L_j$ by $M$. Since $M$ is 
$a$-invariant, 4.12 shows that $M=M^+\oplus M^-$, where
$M^+=L^+\cap M$ and $M^-=L^-\cap M$. It is easy to check (using 4.6 ) that for
any $m\in M^+$ there exists $d\in C_{D_j}(a)$ such that $dD_{j+1}=m$. By the
hypothesis $d$ is of order dividing $n$. Therefore, by Lemma 2.6, $m$ is 
ad-nilpotent of index at most $n$. Similarly one concludes that any element in 
$M^-$ is ad-nilpotent of index at most $n$. 

Since $G=\langle g,h\rangle$, it follows that the subalgebra $L_p(G)$
of $L$ is generated by ${\tilde g}=gD_2$ and ${\tilde h}=hD_2$. The assumption
that $g^a=g^{-1}$ and $h^a=h$ implies that ${\tilde g}\in L^-$, ${\tilde h}\in 
L^+$. By 4.13 we conclude that any commutator in ${\tilde g}, {\tilde h}$ lies
either in $L_j\cap L^+$ or in $L_j\cap L^-$ for some $j$. Combining this with
the observations made in the preceding paragraph, we arrive at the 
conclusion that
any commutator in ${\tilde g}, {\tilde h}$ is ad-nilpotent of index at most $n$.

Arguing like in the proof of Theorem 4.3, remark that the identity 2.9 is
satisfied in $C_L(a)$. Therefore, by Theorem 2.3, $L$ satisfies a certain
polynomial identity depending only on $n$. Theorem 2.1 now shows that
$L_p(G)$ is of $n$-bounded nilpotency class.

By Lemma 4.5 we derive that $G$ contains a characteristic powerful subgroup
$G_1$ of $n$-bounded index. Lemma 4.12 implies that $G_1$ is generated by 
elements of order dividing $n$. Therefore the exponent of $G_1$ divides $n$.
The theorem follows.
$\Box$
\bigskip

It is not clear whether the above theorem can be extended to the case of automorphism of any order prime to that of $G$.
\bigskip

\noindent{\bf Problem 3.}
Let $n$ be a positive integer, $G$ a finite group admitting an automorphism $a$,
of order prime to $|G|$, such that $C_G(a)$ is of exponent dividing $n$.
Suppose that for any $x\in G$ the element $[x,a]=x^{-1}x^a$ has order dividing 
$n$. Is then the exponent of $G$ bounded in terms of $n$ and $|a|$?
\bigskip

\bigskip

\bigskip

\noindent{\bf 5. On centralizers in periodic residually finite groups}
\bigskip

In this section we find some sufficient conditions for a
 periodic residually finite group 
to be locally finite.
As is attested by the groups constructed in \cite{al}, \cite{gol}, \cite{gri}, 
\cite{gs}, \cite{susch}, in general a periodic residually finite group need not
be locally finite. Our theme will be the following.

Given a periodic residually finite group $G$ acted on by a finite group $A$,
under what conditions on $A$ and $C_G(A)$ does it follow that $G$ is locally
finite? Since any subgroup of $G$ acts on $G$ by inner automorphisms, this 
problem includes problems on centralizers of finite subgroups of $G$.

In 1972 V. Shunkov proved that if a periodic group $G$ admits an involutory
automorphism  $a$ with finite centralizer $C_G(a)$ then $G$ contains a solvable
subgroup of finite index \cite{shun}. This was strengthened later by B. Hartley 
and Th. Meixner who showed that $G$ has a nilpotent subgroup of index depending 
only on $|C_G(a)|$ and of nilpotency class at most two \cite{hame}.  Locally 
finite groups
$G$ having an automorphism $a$ of arbitrary prime order $p$ with finite
centralizer $C_G(a)$ have been studied
intensively in seventies and eighties (see for example \cite{ha}). Khukhro
showed that these groups have a nilpotent subgroup of finite index depending
only on $|C_G(a)|$ and on $p$, and nilpotency class depending only on $p$
\cite{khu}.  In general, a very interesting direction in locally finite group 
theory is to classify in some sense locally finite groups $G$ having a finite 
subgroup $A$ such that $C_G(A)$ possesses some prescribed property,
as for example the property to be a linear group \cite{ha2}.

An immediate corollary of the theorem of Shunkov is that $G$ is locally finite.
This part of the theorem, and in fact most difficult part, has no analogue
for periodic groups admitting an automorphism of odd order.
According to G. Deryabina and A. Ol'shanskii, for any
positive integer $n$ which has at least one odd divisor there exists an 
infinite group $G$ having a non-central element of order $n$ such that all
proper subgroups of $G$ are finite \cite{dolsh}. 
Therefore the result of Shunkov cannot be extended to periodic 
groups having an automorphism whose order is not a 2-power.
Moreover, Obraztsov and Miller constructed 
for any (not necessarily distinct) odd primes $p$ and $q$ 
a finitely generated infinite residually finite periodic 
$p$-group admitting a fixed-point-free automorphism of order $q$ \cite{ob}.

On the other hand, in the recent years new means to treat the problem in the 
residually finite case have been found. N. R. Rocco and the author proved in
\cite{rsh} that if a periodic residually finite group $G$ admits an 
automorphism $a$ of
order $2^s$ with $C_G(a)$ finite then $G$ is locally finite. Somewhat later
(\cite{shu3}, \cite{shu4})
the author proved local finiteness of a periodic residually finite group $G$
in the following cases:

1) $G$ admits a 4-group $A$ of automorphisms with $C_G(A)$ finite; or

2) $G$ has no elements of order two and admits an involutory automorphism $a$
such that $C_G(a)$ is abelian.

Then A. Shalev obtained in \cite{sha4} some very general results on local 
finiteness of periodic residually finite groups acted on by a finite 2-group
(see 5.2 and 5.3 bellow). These results were extended in \cite{shu2} to the case
when the acting group is not necessarily of 2-power order.

If $A$ is any finite group, let $q(A)$ denote the maximal prime divisor of 
$|A|$. One of the results described in this section is the following
theorem.
\bigskip

\noindent{\bf Theorem 5.1.} 
Let $G$ be a residually finite group acted on by a finite solvable
group $A$ with $q=q(A)$. Assume that $G$ has no $|A|$-torsion and $C_G(A)$ 
is either solvable or of finite exponent. If any $q-1$ elements of $G$ 
generate a finite solvable subgroup then $G$ is locally finite.
\bigskip

We should mention that for any integer $d\geq 2$
there exist infinite $d$-generator residually finite groups in which all 
$(d-1)$-generator subgroups are finite. The corresponding examples are provided
by Golod's groups \cite{gol}.
 
If under the hypothesis of 5.1 $A$ is a 2-group then the condition imposed on $G$
is that $G$ is merely periodic. This important special case is due to A. Shalev
\cite{sha4}.
\bigskip

\noindent{\bf Theorem 5.2} (Shalev).
Let $G$ be a periodic residually finite group with no 2-torsion acted on by a 
finite 2-group $A$. Suppose the centralizer
$C_G(A)$ is solvable, or of finite exponent. Then $G$ is locally finite.
\bigskip

Other results to be described in this section are as follows.
\bigskip

\noindent{\bf Theorem 5.3} (Shalev).
Let $G$ be a periodic residually finite group having a finite 2-subgroup $A$ 
such that the centralizer $C_G(A)$ is finite. Then $G$ is locally finite.  
\bigskip

\noindent{\bf Theorem 5.4.} 
Let $q$ be a prime, $G$ a residually finite group in which each 2-generator 
subgroup is finite. Suppose that $G$ has a finite $q$-subgroup $A$ 
such that the centralizer $C_G(A)$ is finite. Then $G$ is locally finite.  
\bigskip

Recall that the Fitting height of a finite solvable group $G$ is defined as
the least number $h=h(G)$ such that $G$ possesses a normal series
$$G=G_1\geq G_2\geq\dots\geq G_{h+1}=1,$$ all of whose quotients $G_i/G_{i+1}$
are nilpotent.  Thus, $G$ is nilpotent if and only if $h(G)=1$. In this section 
we will use some deep results on Fitting height of finite solvable groups. 
\bigskip

\noindent{\bf Theorem 5.5} (Thompson, \cite{tho2}). Let $G$ and $A$ be finite
solvable groups such that $(|G|,|A|)=1$. Assume that $A$ acts on $G$ in such a 
manner that $h(C_G(A))=h$. Assume further that the order of $A$ is a product
of $k$ not necessarily distinct primes. Then $h(G)$ is $\{h,k\}$-bounded.
\bigskip

The following lemma
is a corollary of the famous theorem of P. Hall and G. Higman on $p$-length
of a $p$-solvable finite group of given exponent \cite{hahi}.
\bigskip

\noindent{\bf Lemma 5.6} (Shalev, \cite{sha4}). The Fitting height of a finite
solvable group of exponent $n$ is $n$-bounded.
\bigskip

\noindent{\bf Lemma 5.7.} Let $h$ be a positive integer, $G$ a residually 
finite-solvable group such that $h(Q)\leq h$ for any finite quotient
$Q$ of $G$. Then $G$ possesses a normal series
$$G=G_1\geq G_2\geq\dots\geq G_{h+1}=1,$$ all of whose quotients are residually
finite-nilpotent.
\medskip

\noindent{\bf Proof.} Let us use the induction on $h$, the case $h=1$ being 
obvious. Assume that $h\geq 2$, and let $H$ be the intersection of all finite
index normal subgroups $N$ of $G$ such that $h(G/N)\leq h-1$. By the induction
hypothesis 
$G$ possesses a normal series
$$G=G_1\geq G_2\geq\dots\geq G_{h}=H,$$ all of whose quotients are residually
finite-nilpotent. Therefore, it suffices to show that $H$ is residually 
nilpotent. Let $x$ be any non-trivial element of $H$. Since $G$ is residually
finite, there exists a normal subgroup $N$ of finite index in $G$ such that
$x\not\in N$. By hypothesis $G$ possesses a normal series
$$G=N_1\geq N_2\geq\dots\geq N_{h+1}=N,$$ all of whose quotients are 
nilpotent. We note that $h(G/N)=h$ for $N$ does not contain $H$.
Therefore $h(G/N_h)=h-1$ and so $H\leq N_h$. Since the quotient $N_h/N$ is
nilpotent, so is $H/H\cap N$. Thus for any non-trivial element $x\in H$ we can
find a normal in $G$ subgroup $N$ such that $x\not\in N$ and $H/N\cap H$ is
finite and nilpotent. This means that $H$ is residually finite-nilpotent,
as required.
$\Box$
\bigskip

Let $G$ be a periodic group acted on by a finite group $A$. Suppose that $G$
has no non-trivial elements of order dividing that of $A$, and let $N$ be an 
$A$-invariant normal subgroup of $G$. Lemma 4.6 says that if $G$ is finite
then $C_{G/N}(A)=C_G(A)N/N$.
We saw in the previous section that this fact is of fundamental importance
for using Lie ring methods in the study of finite groups $G$ having automorphisms
of coprime order. The applicability of Lie ring methods to the study of infinite
periodic groups $G$ acted on by a finite group $A$ depends ultimately on how
successful one is in extending Lemma 4.6 to the case when $G$ is allowed to be
infinite periodic. Since in general
$C_{G/N}(A)\neq C_G(A)N/N$, one has to impose additional conditions
on $G$ and $A$.

Given a positive integer $n$, a group $G$ is said to be $n$-finite if any
$n$-generator subgroup of $G$ is finite. Thus, a group is 1-finite if and only
if it is periodic. It was proved in \cite{shu1} that the equality
$C_{G/N}(A)=C_G(A)N/N$ holds whenever $A$ is a 2-group. More generally, we have.
\bigskip

\noindent {\bf Lemma 5.8.} Let $A$ be a finite solvable group with $q=q(A)$
acting on a $(q-1)$-finite group $G$ with no $|A|$-torsion. Let $N$ be a normal 
$A$-invariant subgroup of $G$. Then $C_{G/N}(A)=C_G(A)N/N$.
\medskip

\noindent{\bf Proof.} It suffices to show that any $A$-invariant coset $xN$
contains an element from $C_G(A)$. Suppose first that $A$ is of order $q$ and
let $a$ be a generator of $A$. Set 
\begin{center}
$x_0=x^{-1}x^a$, $x_1=x_0^a$, $\dots$, $x_{q-1}=x_{q-2}^a$.
\end{center}
Then all 
$x_0,x_1,\dots, x_{q-1}$ 
lie in $N$ and $x_0x_1\dots x_{q-1}=1$. It follows 
that the subgroup
$F=\langle x_0,x_1,\dots ,x_{q-1}\rangle$ 
is generated by at most $q-1$ elements. Hence
$F$ is finite. We note that $F$ is $A$-invariant and consider the natural 
split extension $FA$.

Since $G$ has no $|A|$-torsion, it is easy to see that $A$ and 
$B=\langle ax_0^{-1} \rangle$ are Sylow $q$-subgroups of $FA$. Hence there 
exists an element $y\in F$ such that $B=A^y$. Therefore $xy^{-1}\in N_G(A)=
C_G(A)$ and $x\in C_G(A)N$, as required.

Now let $A$ be of non-prime order. Let $D$ be a non-trivial proper normal 
subgroup of $A$. Since $D$ is of order less than that of $A$, we can assume by 
induction that $xN\cap C_G(D)\neq \emptyset$. Therefore without any loss of
generality we can assume $x\in C_G(D)$. Let $H=C_G(D)$, ${\bar A}=A/C_A(H)$,
$M=N\cap H$. The coset $xM$ is obviously ${\bar A}$-invariant. Arguing by 
induction on $|A|$ and using that ${\bar A}$ has order less than that of $A$
we can assume that $xM$ contains an element $z\in C_H({\bar A})$. Now it remains
to notice that $zN=xN$ and $z\in C_G(A)$.
$\Box$
\bigskip

\noindent {\bf Lemma 5.9. } Let $q$ be a prime, $A$ a finite $q$-group acting 
on a residually finite group $G$.  Suppose that $G$ is $(q-1)$-finite and
$C_G(A)$ has no $q$-torsion. Then $G$ has no $q$-torsion.
\medskip

\noindent {\bf Proof.}  Suppose that the lemma is false and assume first that 
$A$ is of order $q$. Let $a$ be a generator of $A$. Since $G$ is residually 
finite, we can choose a normal $A$-invariant subgroup $N$ such that $G/N$ is 
a finite group whose order is divisible by $q$. It follows that $A$ 
centralizes some element $xN$ of order $q$ in $G/N$.  Then $x^{-1}x^a\in N$.
Set
\begin{center}
$x_0=x^{-1}x^a$, $x_1=x_0^a$, $\dots$, $x_{q-1}=x_{q-2}^a$.
\end{center}
Arguing like in the previous lemma we observe that 
$F=\langle x_0,x_1,\dots ,x_{q-1}\rangle$ is finite. If $q$ divides $|F|$ 
then obviously $A$ must centralize some elements of order $q$ in $F$. This 
yields a contradiction, so assume that $|F|$ is prime to $q$. Then 
$\langle a{x_0}^{-1}\rangle$ is a Sylow $q$-subgroup of $FA$. Therefore it is 
a conjugate of $A$. Hence there exists an element $y\in F$ such that 
$y^{-1}ay=x^{-1}ax$. Then $z=xy^{-1}\in C_G(A)$. Since the image of 
$z$ in $G/N$ is of order $q$, it follows that the order of $z$ is divided by 
$q$. Therefore $C_G(A)$ contains an element of order $q$, a contradiction.

Suppose now that $A$ is of order $q^n$ and use induction on $n$. Let $a$ be an
element of prime order in $Z(A)$. Set $C=C_G(a)$. If $C$ is $q$-torsion free
then by the previous paragraph so is $G$. Assume that $C$ has non-trivial 
$q$-elements. Since $A$ induces an automorphism group of $C$ whose order is 
strictly less than that of $A$, the induction hypothesis implies that some of 
$q$-elements of $C$ must lie in $C_G(A)$.
$\Box$
\bigskip

Obviously the conclusions of Lemmas 5.8 and 5.9 remain true if we replace the 
assumption that $G$ is $(q-1)$-finite by the assumption that the semidirect 
product $GA$ is 2-finite. Really, keeping notation like in Lemma 5.8 let us 
note that if $GA$ is 2-finite then $F$ is finite because 
$F\leq\langle x,a\rangle$. Leaving 
other parts of proofs unchanged we reach the following results. 
\bigskip

\noindent {\bf Lemma 5.8$'$.} Let $A$ be a finite solvable group 
acting on a group $G$ with no $|A|$-torsion. Suppose that $GA$ is 2-finite.
Let $N$ be a normal $A$-invariant subgroup of $G$. Then $C_{G/N}(A)=C_G(A)N/N$.
\bigskip

\noindent {\bf Lemma 5.9$'$.} Let $q$ be a prime, $A$ a finite $q$-group acting 
on a residually finite group $G$.  Suppose that $GA$ is 2-finite and
$C_G(A)$ has no $q$-torsion. Then $G$ has no $q$-torsion.
\bigskip

\noindent{\bf Proposition 5.10.} Let $p$ be a prime and $G$ a periodic
residually finite $p$-group acted on by a finite solvable group $A$ whose order 
is prime to $p$. Suppose that $q=q(A)$ and $G$ is $(q-1)$-finite. If $C_G(A)$ 
satisfies a non-trivial identity then $G$ is locally finite.
\smallskip

\noindent{\bf Proof.} Since any finite set of elements of $G$ is contained
in a finitely generated $A$-invariant subgroup, we can assume that $G$ is
finitely generated. Let $L=L_p(G)$ be the Lie algebra associated with the
$p$-dimension central series of $G$. We regard $A$ as a group acting on $L$.
Just as in Section 4 (but using Lemma 5.8 in place of Lemma 4.6) we remark
that $C_L(A)=L_p(G,C_G(A))$ and so, because $C_G(A)$ satisfies a non-trivial
identity, by Proposition 2.7 $C_L(A)$ is PI. Since $G$ is periodic, it follows
from Lemma 2.6 that $L$ has a finite set of generators in which every
commutator is ad-nilpotent. Consequently, by Theorem 1.5
combined with Theorem 2.2 $L$ is nilpotent.

Let $g_1,g_2,\dots,g_m$ be some generating set of $G$, let $c$ be the nilpotency
class of $L$ and write $\rho_1,\rho_2,\dots,\rho_s$ for the list of all simple
commutators of weight at most $c$ in $g_1,g_2,\dots,g_m$. Choose $K$ to be the 
maximum of orders of $\rho_i$. Then, by 2.14, $D_j(G)$ has index at most $K^s$
for any $j$. Since $G$ is residually-$p$, it follows that the order of $G$ is
at most $K^s$. The proof is complete.
$\Box$
\bigskip

\noindent{\bf Proof of Theorems 5.1 and 5.2.}
Let us note that $G$ is residually solvable. Indeed, if $q=2$ then $G$ has no
2-torsion and so $G$ is residually solvable by the Feit-Thompson Theorem 
\cite{feth}. If $q\geq 3$ then any two elements of $G$ generate a finite solvable
subgroup. Since any simple finite group can be generated by two elements, it 
follows that $G$ is residually solvable. 

Assume that $C_{G}(A)$ is of finite 
exponent. Let $N$ be any $A$-invariant normal subgroup of finite index in
$G$ and $Q={G}/N$. Then $Q$ is solvable because $G$ is residually so.
By Lemma 5.8 $A$ acts on $Q$ in such a way that $C_Q(A)$ is of exponent at most
that of $C_{G}(A)$. By Lemma 5.6 the 
Fitting height $h(C_Q(A))$ of $C_Q(A)$ is bounded in terms of the exponent of 
$C_Q(A)$. Theorem 5.5 now shows that the Fitting height 
$h(Q)$ of $Q$ is bounded in terms of the exponent of $C_{G}(A)$ and $|A|$. 
Thus, $h(Q)$ is bounded by a number $h$ which does 
not depend on $Q$. Therefore, by Lemma 5.7, $G$ possesses a normal $A$-invariant
series of length at most $h+1$ all of whose quotients are residually nilpotent.
Induction on $h$ shows that without any loss of generality $G$ can be assumed 
residually nilpotent. Then $G$ is a direct product of $A$-invariant 
$p$-subgroups and finiteness of $G$ follows from Proposition 5.10. 

If $C_{G}(A)$ is solvable and $Q$ has the same meaning as above then
$h(C_Q(A))$ is at most the derived length of $C_{G}(A)$ and we can repeat
the argument. This completes the proof.
$\Box$
\bigskip

Using Lemma 5.9 and 5.9$'$ we can in some situations derive local finiteness of 
$G$ even without requiring that $G$ has no $|A|$-torsion.
\bigskip

\noindent{\bf Proof of Theorems 5.3 and 5.4.} Since $G$ is residually finite, 
we can choose a normal subgroup $H$ of finite index in $G$ such that
$H\cap C_G(A)=1$. It suffices to prove that $H$ is locally finite. Lemmas 5.9 
and 5.9$'$ imply that $H$ has no $|A|$-torsion. Hence, by 5.8 (or 5.8$'$), $A$
acts fixed-point-freely on every $A$-invariant finite quotient $Q$ of $H$.
A well-known corollary of the classification of simple finite groups says
that any finite group admitting a fixed-point-free automorphism group of
coprime order is solvable. We now conclude that $Q$ is necessarily solvable.
Hence $H$ is residually solvable. This places us in a position to apply
Theorem 5.1 (or Theorem 5.2 in case $q$=2) and derive that $G$ is locally 
finite.
$\Box$
\bigskip

Since there is no restriction on the identity satisfied by $C_G(A)$ in 
Proposition 5.10, it is not unreasonable to conjecture that the following
problem can be answered positively.
\bigskip

\noindent{\bf Problem 4.} Let $A$ be a finite 2-group 
acting on a periodic residually finite group $G$ which has no 2-torsion.
Assume that $C_G(A)$ satisfies a non-trivial identity. Does it follow that $G$
is locally finite?
\bigskip

\vskip3ex

{\small


\begin{thebibliography}{12} 
\bibitem{al} S. V. Aleshin, Finite automata and the Burnside problem for 
periodic groups, Math. Notes {\bf 11} (1972), 199--203.

\bibitem{bz} Yu. A. Bahturin, M. V. Zaicev, 
Identities of graded 
algebras, {\it to appear in J. Algebra}. 
\bibitem{dolsh} G. S. Deryabina and A. Yu. Ol'shanskii, Subgroups of 
quasifinite groups, Uspekhi Math. Nauk, {\bf 41} (1986), 169--170
(in Russian).
\bibitem{ddms}  J. D. Dixon, 
M. P. F. du\,\,Sautoy, A. Mann, 
D. Segal, ``Analytic pro-$p$ groups" (London Math. 
Soc. 
Lecture Note Series {\bf 157}), Cambridge Univ. Press., 
1991. 
\bibitem{feth} W. Feit, J. G. Thompson, Solvability of 
groups  of odd order, Pacific J. Math. {\bf 13} (1963), 773-1029. 
\bibitem{gol} E. S. Golod, On nil-algebras and residually finite groups,
Izvestia Akad. Nauk SSSR, Ser Mat., {\bf 28} (1964), 273--276

\bibitem{go} D. Gorenstein, ``Finite Groups", 
Harper and 
Row, New 
York, 1968. 
\bibitem{gri} R. I. Grigorchuk, On the Burnside problem for periodic
groups, Funct. Anal. Appl., {\bf 14} (1980), 53--54
\bibitem{gr} K. Gruenberg, The Engel elements of a soluble group, Illinois
J. Math. {\bf 3} (1959), 151-168.
\bibitem{gs} N. Gupta and S. Sidki, On the Burnside problem for
periodic groups, Math. Z., {\bf 182} (1983), 385--386.
\bibitem{hahi} P. Hall, G. Higman, The $p$-length of a $p$-soluble group
and reduction theorems for Burnside's problem, Proc. London. Math. Soc. (3)
{\bf 6} (1956), 1--42.

\bibitem{ha4} B. Hartley, Subgroups of finte index in
profinite groups, Math. Z., {\bf 168} (1979), 71--76. 
\bibitem{ha} B. Hartely, Centralizers in locally finite groups,
in ``Proc. 1st Bressanone Group Theory Conference 1986", Lecture
Notes in Math. 1281, Springer, 1987, pp. 36--51.
\bibitem{ha2} B. Hartley, Simple locally finite groups, in ``Finite and Locally Finite 
Groups", NATO ASI Series, {\bf 471}, Kluwer Academic Publishers, 1995, 1--45.
\bibitem{hame} B. Hartley, T. Meixner, Periodic groups in which the 
centralizer of an involution has bounded order,  Arch. Math. (Basel), {\bf 
 36} (1981), 211--213. 
\bibitem{hi} G.\,Higman,  Lie ring methods in the 
theory of finite nilpotent groups, in ``Proc. Intern. Congr. Math. 
Edinburgh, 1958", Cambridge Univ. Press, 1960, 307--312. 
\bibitem{hu} B. Huppert, ``Endliche Gruppen I", Springer Verlag, Berlin, 1967. 
\bibitem{hb} B. Huppert, N. Blackburn, ``Finite Groups II", Springer Verlag, Berlin, 1982. 
\bibitem{khu} E. I. Khukhro, Groups and Lie rings admitting an almost regular 
automorphism of prime order, Math. USSR Sbornik, {\bf 71} (1992), 51--63. 
\bibitem{khu1} E. I. Khukhro,  ``Nilpotent Groups and their Automorphisms", 
de\,\,Gruyter--Verlag, Berlin, 1993.
\bibitem{khsh} E. I. Khukhro and P. V. Shumyatsky, Bounding the exponent of a
finite group with automorphisms, Preprint, {\it to appear in J. Algebra}.
\bibitem{ko} 
A. I. Kostrikin, On the Burnside problem, Izv. AN SSSR, Ser. Mat. 
{\bf 23} (1959), 3--34. 
\bibitem{la} M. Lazard,  Sur les groupes nilpotents et les anneaux de Lie, 
Ann. Sci. \'Ecole Norm. Supr. {\bf 71} (1954), 101-190.
\bibitem{l2} M.\,Lazard, Groupes analytiques $p$-adiques, Publ. Math. Inst. 
Hautes \'Etudes Sci., {\bf 26} (1965), 389--603. 
\bibitem{li}  V. Linchenko, Identities of Lie algebras with actions of 
Hopf algebras, {\it to appear in Commun. Algebra.} 
\bibitem{lbmn} A. Lubotzky, A. Mann, Powerful $p$-groups. I: finite groups, 
J. Algebra, {\bf 105} (1987), 484--505; II: $p$-adic analytic groups, {\it ibid.}, 506--515. 
\bibitem{ma} W. Magnus, A connection between the 
Baker--Hausdorff formula and a problem of Burnside, Ann. of Math. (2), 
{\bf 52} (1950), 111--126.
\bibitem{ob}  C. F. Miller III and V. N. Obraztsov, Infinite periodic residually
finite groups with all finite subgroups cyclic, in preparation.
\bibitem{ol} A. Yu. Ol'shanskii, ``Geometry of defining relations in groups",
Math. Appl. (Soviet Ser.) {\bf 70}(1991).
\bibitem{rsh} N. Rocco, P. Shumyatsky, On periodic groups having almost regular 
$2$-ele\-ments,  Proc. Edinburgh Math. Soc., {\bf 41}(1998), 385--391.
\bibitem{sa}  I. N. Sanov, Establishment of a 
connection between periodic groups 
with period a prime number and Lie rings, Izv. Akad. Nauk SSSR 
Ser. Mat., {\bf 16} (1952), 23--58 (Russian). 
\bibitem{sha1} A. Shalev, Characterization of $p$-adic analytic groups in terms 
of wreath products, J. Algebra, {\bf 145}(1992), 204--208.
\bibitem{sha5}  A. Shalev, On almost fixed point free 
automorphisms,  J. Algebra, {\bf 157} (1993), 271--282. 
\bibitem{sha6} A. Shalev, Automorphisms of finite groups of bounded rank,
Israel J. Math., {\bf 82}(1993), 395--404.
\bibitem{sha2} A. Shalev, Finite $p$-groups, in ``Finite and Locally Finite 
Groups", NATO ASI Series, {\bf 471}, Kluwer Academic Publishers, 1995, 401--450.
\bibitem{sha4} A. Shalev, Centralizers in residually 
finite torsion groups, {\it to appear in Proc. Amer. Math. Soc.} 
\bibitem{shu1} P. V. Shumyatsky,
Groups with regular elementary 2-groups of automorphisms, Algebra and Logic,
{\bf 27} (1988), 447--457.
\bibitem{shu0} P. Shumyatsky, On periodic solvable groups
having automorphisms with
nilpotent fixed-point groups, Israel J. Math., {\bf 87} (1994), 111--116.
\bibitem{shu3} P. Shumyatsky, On groups having a 
four-subgroup with finite centralizer, {\it to appear in Quart. J. Math.} 
\bibitem{shu4}  P. Shumyatsky, Nilpotency of some Lie algebras associated with 
$p$-groups, {\it Preprint of the Univ. of Brasilia}, 1997. 
\bibitem{shu2} P. Shumyatsky, Centralizers in groups with finiteness 
conditions, {\it to appear in J. Group Theory}. 
\bibitem{shu0} P. Shumyatsky, On groups with commutators of
bounded order, {\it to appear in Proc. Amer. Math. Soc}. 
\bibitem{shun} V.P. Shunkov, On periodic groups with an almost regular
involution, Algebra and Logic, {\bf 11} (1972), 260--272.
\bibitem{susch} V.I. Sushchansky, Periodic $p$-elements of permutations
and the general Burnside problem, Dokl. Akad. Nauk SSSR, {\bf 247} (1979),
447--461.
\bibitem{tho2} J.G. Thompson, Automorphisms of solvable groups,
J. Algebra, {\bf 1} (1964), 259-267.
\bibitem{tits} J. Tits, Free subgroups in linear groups, J. Algebra, {\bf 20}
(1972), 250-270.
\bibitem{vl} 
M. R. Vaughan-Lee, ``The Restricted Burnside Problem", 2-nd edition, Oxford
University Press, Oxford, 1993.
\bibitem{wa} J. N. Ward, Automorphisms of finite groups and their fixed-point
groups, J. Austral. Math. Soc., {\bf 9} (1969), 467--477.
\bibitem{w} J.S. Wilson, Two-generator conditions for residually finite
groups, Bull. London Math. Soc., {\bf 23} (1991), 239-248. 
\bibitem{wize} J.S. Wilson and E. Zelmanov, Identities for Lie algebras
of pro-p groups, J. Pure Appl. Algebra, {\bf 81} (1992), 103-109. 
\bibitem{za} H. Zassenhaus, Ein Verfahren, jeder endlischen Gruppe einen
Lie-Ring mit der Charakteristiki p zuzuordnen,
Abh. Math. Seminar Hans. Univ. Hamburg {\bf 13}(1940), 200--207
\bibitem{ze1} E. Zelmanov,  The solution of the restricted Burnside problem
for groups of odd exponent, Math. USSR Izv. {\bf 36} (1991), 41-60.
\bibitem{ze2} E. Zelmanov,  The solution of the restricted Burnside problem
for 2-groups , Math. Sb. {\bf 182} (1991), 568-592.
\bibitem{zelm} 
E. Zelmanov,  ``Nil Rings and Periodic Groups", The Korean Math. Soc. Lecture
Notes in Math., Seoul, 1992.
\bibitem{ze3} E. I. Zelmanov, Lie ring methods in the 
theory of nilpotent groups,  in ``Proc. Groups'93/St. Andrews, vol. 2" 
(London Math. Soc. Lecture Note Ser. {\bf 212}),  
Cambridge Univ. Press, 1995, 567--585. 




\end{thebibliography}
\end{document}